\documentclass[a4paper,12pt]{article}
\usepackage[english]{babel}
\usepackage{amsfonts,amssymb,amsmath}

\begin{document}

 \centerline{\bf\Large Long-time asymptotics of solutions} 
 \vskip 3mm 
 \centerline{\bf\Large of the heat equation} 

 \vskip 7mm 
 \centerline{\bf S.V.~Zakharov}

\begin{center}
Institute of Mathematics and Mechanics,\\
Ural Branch of the Russian Academy of Sciences,\\
16, S.\,Kovalevskaja street, 620990, Ekaterinburg, Russia
\end{center}

\vskip 5mm 

\textbf{Abstract.}
The long-time asymptotics of solutions of the Cauchy problem
for the heat equation are constructed 
in the case when the initial function 
at infinity has power asymptotics. 

\vskip 3mm 

Key words: heat equation, Cauchy problem, long-time asymptotics, Hermite functions.

Mathematics Subject Classification: 35K05, 35K15.
 
\vskip 7mm
\section{ One-dimensional heat equation}
 
 Consider the Cauchy problem for
the one-dimensional heat equation
with a bounded and continuous
 initial function $\Lambda:\mathbb{R} \to \mathbb{R}$: 
 \begin{eqnarray} 
 \label{eq1} 
 \frac{\partial u}{\partial t} = 
 \frac{\partial^2 u}{\partial x^2},& 
 \quad t\geqslant 0, \\ 
 \label{ic1}\phantom{\frac{1}{1}} 
 u(x,0) = \Lambda (x),& 
 \quad x\in\mathbb{R}. 
 \end{eqnarray} 
 
 The solution of this problems is written in the form
of convolution 
\begin{equation}\label{us1} 
 u(x,t) = \frac{1}{2\sqrt{\pi t}} 
 \int\limits_{-\infty}^{+\infty} 
 \Lambda(s) \exp\left\{ - \frac{(s-x)^2}{4t} \right\}\,
ds. 
 \end{equation} 
 
Investigation of the asymptotic behavior of
integral~(\ref{us1}) in addition to direct applications 
to physical processes of heat conduction and diffusion
is of independent interest for asymptotic analysis,
because the necessity to solve such problems 
arises under applying the matching method~\cite{ib}. 
 
We assume that there hold the following asymptotic
relations: 
 \begin{equation}\label{le1} 
 \Lambda(x) = \sum\limits_{n=0}^{\infty} 
 \frac{\Lambda^{\pm}_n}{x^{n}}, 
 \qquad 
 x\to \pm\infty. 
 \end{equation} 
 
 Let us represent function~(\ref{us1}) in the
form $$ 
 u(x,t)=U^-_1(x,t) + U^-_0(x,t) + U^+_0(x,t) + U^+_1(x,t), 
 $$ 
 where $$ 
 U^-_1(x,t) = \int\limits_{-\infty}^{-\sigma} \dots \,
ds, 
 \qquad 
 U^-_0(x,t) = \int\limits_{-\sigma}^{0} \dots \,
ds, 
 $$ 
 $$ 
 U^+_0(x,t) = \int\limits_{0}^{\sigma} \dots \,
ds, 
 \qquad 
 U^+_1(x,t) = \int\limits_{\sigma}^{+\infty} \dots \,
ds, 
 $$ 
 $$ 
 \sigma = (x^2 + t)^{p/2}, 
 \qquad 
 0 < p < 1. 
 $$ 
 
 In the integral $U^+_1(x,t)$
 we make the change $s = 2z\sqrt{t}$. 
 Setting $$ 
 \mu= \frac{\sigma}{2\sqrt{t}}, 
 \qquad 
 \eta = \frac{x}{2\sqrt{t}} 
 $$ 
 and using condition~(\ref{le1}), we obtain 
$$ 
 U^+_1(x,t) = \frac{1}{\sqrt{\pi}} 
 \int\limits_{\mu}^{+\infty} 
 \Lambda(2z\sqrt{t}) e^{-(z-\eta)^2}
dz = 
 $$ 
 \begin{equation}\label{up1i} 
 = \sum\limits_{n=0}^{N-1} 
 \frac{\Lambda^+_n}{\sqrt{\pi} 2^n} 
 t^{-n/2} 
 \int\limits_{\mu}^{+\infty} z^{-n} 
 e^{-(z-\eta)^2} dz + 
 \int\limits_{\mu}^{+\infty} 
 R_N(z\sqrt{t}) e^{-(z-\eta)^2} dz, 
 \end{equation} 
 where $$ 
 |R_N(s)|\leqslant K_N s^{-N}. 
 $$ 
 From the last inequality we obtain
the estimate
 \begin{equation}\label{ern} 
 \left| \int\limits_{\mu}^{+\infty} 
 R_N(z\sqrt{t}) e^{-(z-\eta)^2} dz 
 \right| 
 \leqslant k_N t^{-N/2} 
 \int\limits_{\mu}^{+\infty} z^{-N}
dz = 
 \tilde{k}_N \mu \sigma^{-N}. 
 \end{equation} 
 
Let us extract the dependence on the parameter $\mu$ in
the integral $$ 
 \int\limits_{\mu}^{+\infty} z^{-n} 
 e^{-(z-\eta)^2} dz. 
 $$ 
 
 For $n=0$ and $t\geqslant |x|^{\alpha}$, 
 $1 + p < \alpha < 2$, we have $$ 
 \int\limits_{\mu}^{+\infty} 
 e^{-(z-\eta)^2} dz = 
 \int\limits_{0}^{+\infty} 
 e^{-(z-\eta)^2} dz - \int\limits_{0}^{\mu} 
 e^{-(z-\eta)^2} dz = 
 $$ 
 \begin{equation}\label{ei0} 
 =\sqrt{\pi}\, \mathrm{erfc}(-\eta) + 
 \sum\limits_{r=1}^{N-1} \mu^r 
 e^{-\eta^2} P_{r-1}(\eta) + 
 O(\sigma^{-\gamma N}), 
 \qquad 
 \sigma\to \infty, 
 \end{equation} 
 where $P_m(\eta)$ denotes polynomials of the~$m$th
degree in $\eta$, 
 $$ 
 \gamma = \frac{\alpha}{2p} - 1 > 0, 
 \qquad 
 \mathrm{erfc} (x) = \frac{1}{\sqrt{\pi}} 
 \int\limits_{x}^{+\infty} e^{-s^2}
ds. 
 $$ 
To obtain (\ref{ei0}) we make use of 
 estimates 
\begin{equation}\label{ms} 
 \mu = O\left(t^{\frac{2p-\alpha}{2\alpha}} \right), 
 \qquad 
 \sigma = O\left(t^{p/\alpha} \right), 
 \qquad 
 \mu = O\left( \sigma^{-\gamma} \right)
 \end{equation} 
as $\sigma\to \infty$ on the set $$ 
 T_{\alpha} = \{ (x,t)\, : \, x\in\mathbb{R},\
t \geqslant |x|^{\alpha} \} .
 $$ 

 However, it is seen that~(\ref{ei0}) remains valid
 on the set $$ 
 X_{\alpha} = \{ (x,t)\, : \, x\in\mathbb{R},\ 0 <
t < |x|^{\alpha} \}
 $$ 
as well, where 
$$ 
 \mu = o(\eta), 
 \qquad 
 |\eta| \geqslant \frac{|x|^{1-\frac{\alpha}{2}}}{2} 
 \to \infty 
 \qquad \mbox{for } 
 \quad \sigma\to \infty. 
 $$ 
 
 For the calculation of $U^-_1(x,t)$ 
it is useful an analogous
 relation $$ 
 \int\limits_{-\infty}^{-\mu} 
 e^{-(z-\eta)^2} dz 
 =\sqrt{\pi}\, \mathrm{erfc}(\eta) + 
 \sum\limits_{r=1}^{N-1} \mu^r 
 e^{-\eta^2} P_{r-1}(\eta) + 
 O(\sigma^{-\gamma N}). 
 $$ 
 
 For $n\geqslant 1$ we have $$ 
 \int\limits_{\mu}^{+\infty} z^{-n} 
 e^{-(z-\eta)^2} dz = 
 \int\limits_{1}^{+\infty} z^{-n} 
 e^{-(z-\eta)^2} dz + 
 $$ 
 $$ 
 + \int\limits_{\mu}^{1} \Psi_n(z,\eta)
dz + 
 e^{-\eta^2} \sum\limits_{r=0}^{n-1} P_r(\eta) 
 \int\limits_{\mu}^{1} z^{r-n} dz, 
 $$ 
 where 
\begin{equation}\label{psin} 
 \Psi_n(z,\eta) = z^{-n} \left[ 
 e^{-(z-\eta)^2} - e^{-\eta^2} \sum\limits_{r=0}^{n-1}
z^{r} P_{r}(\eta) 
 \right], 
 \end{equation} 
 and the sum in $r$ is a partial sum of the Taylor  series 
for the function $\exp (2z\eta - z^2)$ in the variable~$z$. 
 Thus, $$ 
 \int\limits_{\mu}^{+\infty} z^{-n} 
 e^{-(z-\eta)^2} dz = 
 \int\limits_{1}^{+\infty} z^{-n} 
 e^{-(z-\eta)^2} dz + 
 e^{-\eta^2} P_{n-1}(\eta) \ln\mu + 
 $$ 
 \begin{equation}\label{ein} 
 +\sum\limits_{r=0}^{n-2} \mu^{r-n+1}
e^{-\eta^2} P_r(\eta) 
 + \int\limits_{0}^{1}\Psi_n(z,\eta)
dz 
 - \int\limits_{0}^{\mu}\Psi_n(z,\eta)
dz. 
 \end{equation} 
 From formula~(\ref{psin}) we conclude
that $\Psi_n(z,\eta)$ has no singularities; 
 from the same formula we obtain $$ 
 \int\limits_{0}^{\mu}\Psi_n(z,\eta)
dz 
 = \sum\limits_{r=1}^{N-1} \mu^{r} e^{-\eta^2}P_{r-1+n}(\eta)+ 
 O(\sigma^{-\gamma N}). 
 $$

 Substituting relation~(\ref{ei0}) and~(\ref{ein}) 
 in formula~(\ref{up1i}) and using the estimate~(\ref{ern}), 
 we obtain $$ 
 U^+_1(x,t) = \Lambda^+_0 \mathrm{erfc}(-\eta) + 
 \phantom{===================} 
 $$ 
 \begin{equation}\label{up1f} 
 +\sum\limits_{n=1}^{N-1} 
 \frac{\Lambda^+_n}{\sqrt{\pi} 2^{n}}
t^{-n/2} 
 \left[ \int\limits_{1}^{+\infty} z^{-n}
e^{-(z-\eta)^2} dz + 
 \int\limits_{0}^{1} \Psi_n(z,\eta)
dz \right] + 
 \end{equation} 
 $$ 
 \phantom{==================} 
 + V_1(\mu,\eta,t) + O(\sigma^{-\gamma N}). 
 $$ 
 Here and further, by $V_i(\mu,\eta,t)$ 
we denote  expressions of the form $$ 
 e^{-\eta^2} \sum\limits_{r^2_s+q^2_s\neq 0} 
 b_s \eta^{m_s} t^{l_s} \mu^{r_s} \ln^{q_s}\mu. 
 $$ 
 
 Since $|x|\leqslant \sigma^{1/p}$, 
 and $t\geqslant \sigma^{\alpha/p}$ on
the set $T_{\alpha}$, 
 we conclude that for $0\leqslant s \leqslant \sigma$ 
 there hold the estimates $$ 
 \frac{xs}{t} = O(\sigma^{-\rho}), 
 \qquad 
 \frac{s^2}{t} = O(\sigma^{-2\rho}), 
 \qquad 
 \rho = \frac{\alpha-1}{p}-1 > 0. 
 $$ 
 Using these estimates, we represent
the integral $U^+_0 (x,t)$ 
 in the following form: 
 $$ 
 U^+_0 (x,t) = \frac{1}{2\sqrt{\pi t}} 
 \int\limits_{0}^{\sigma} \Lambda(s) 
 \exp \left\{ -\frac{(s-x)^2}{4t} \right\}
ds = 
 $$ 
 $$ 
 = \frac{e^{-\eta^2}}{2\sqrt{\pi t}} 
 \int\limits_{0}^{\sigma} \Lambda (s) 
 \sum\limits_{m=0}^{N-1} \frac{1}{m!} 
 \left( \frac{\eta s}{\sqrt{t}} - \frac{s^2}{4t}\right)^m ds + 
 O( \sigma^{- \rho N}). 
 $$ 
Because of the factor $\exp(-\eta^2)$ 
the estimate of the remainder 
is valid on the set $X_{\alpha}$ as well. 
Expanding the power and 
changing the order of summation, we
obtain $$ 
 U^+_0(x,t) 
 = \frac{e^{-\eta^2}}{2\sqrt{\pi t}} 
 \sum\limits_{m=0}^{N-1} t^{-m/2} 
 \sum\limits_{k=0}^{m} 
 \frac{t^{-k/2} \eta^{m-k} }{4^k k! (m-k)!} 
 \int\limits_{0}^{\sigma} s^{m+k} \Lambda(s)
ds+ 
 O( \sigma^{- \rho N}) = 
 $$ 
 
 $$ 
 = \sum\limits_{n=1}^{N-1} t^{-n/2} 
 \sum\limits_{k=0}^{[(n-1)/2]} 
 \frac{\eta^{n-2k-1} e^{-\eta^2}}{2\sqrt{\pi}4^k
k! (n-2k-1)!} 
 \int\limits_{0}^{\sigma} s^{n-1} \Lambda(s)
ds + 
 O(\sigma^{- \rho N}). 
 $$ 
 
 Let us transform the integral as follows: 
 $$ 
 \int\limits_{0}^{\sigma} s^{n-1} \Lambda(s)
ds = 
 \int\limits_{0}^{1} s^{n-1} \Lambda(s)
ds + 
 \int\limits_{1}^{\sigma}\left( 
 \Lambda^+_0 s^{n-1} + \dots + \Lambda^+_{n-1} + 
 \frac{\Lambda^+_n}{s} \right) ds + 
 $$ 
 $$ 
 + \int\limits_{1}^{\sigma} s^{n-1} \left[ 
 \Lambda(s) - \Lambda^+_0 - \dots - \frac{\Lambda^+_{n-1}}{s^{n-1}} - \frac{\Lambda^+_n}{s^{n}} \right] ds = 
 $$ 
 $$ 
 = \int\limits_{0}^{1} s^{n-1} \Lambda(s)
ds + 
 \int\limits_{1}^{+\infty}\Phi^+_n(s)
ds - \sum\limits_{m=1}^{n} \frac{\Lambda^+_{n-m}}{m} + 
 $$ 
 \begin{equation}\label{is1} 
 + \Lambda^+_{n} \ln\sigma + 
 \sum\limits_{m=1}^{n} \frac{\Lambda^+_{n-m}}{m}\sigma^{m} - \int\limits_{\sigma}^{+\infty} \Phi^+_n(s) ds, 
 \end{equation} 
 where $$ 
 \Phi^+_n(s) = s^{n-1} \left[ \Lambda(s) - \sum\limits_{m=0}^{n} \frac{\Lambda^+_m}{s^{m}} \right]. 
 $$ 
 From condition~(\ref{le1}) we obtain the
estimate $$ 
 |\Phi^+_n(s)| \leqslant C_n s^{-2}, 
 \qquad s\geqslant 1; 
 $$ 
 hence, $$ 
 \int\limits_{\sigma}^{+\infty} \Phi^+_n(s)
ds = 
 \sum\limits_{m=1}^{N-1} \varphi_{n,m} \sigma^{-m} + 
 O(\sigma^{-N}), 
 \qquad \sigma\to \infty. 
 $$ 
 Taking into account this relation and
substituting $\sigma=2\mu\sqrt{t}$ in (\ref{is1}), we
obtain $$ 
 \int\limits_{0}^{\sigma} s^{n-1} \Lambda(s)
ds = 
 I_n + \frac{\Lambda^+_n}{2} \ln t + 
 \Lambda^+_n \ln\mu + 
 \sum\limits_{r\neq 0} a_{n,r} \mu^{r}
t^{r/2} 
 + O\left( \sigma^{-N}\right), 
 $$ 
 where $I_n$ are constants. 
 Then 
\begin{equation}\label{up0f} 
 U^+_0(x,t) = \sum\limits_{n=1}^{N-1} 
 t^{-n/2} e^{-\eta^2} 
 [P_{n-1}(\eta) + \widetilde{P}_{n-1}(\eta) \ln
t] 
 +V_{0}(\mu,\eta,t) + O(\sigma^{-\rho N}). 
 \end{equation} 
 
 Expressions $V_{0}(\mu,\eta,t)$ and $V_{1}(\mu,\eta,t)$ 
 in formulas~(\ref{up0f}) and~(\ref{up1f}) 
 asymptotically diminish.
Thus, we arrive at the following statement
(details of the proof see~in~\cite{hd}). 
 
\vskip 5mm 
 \textbf{Theorem~1.} 
 {\it 
 As $|x| + t \to \infty$ the asymptotics of the
solution to problem $(\ref{eq1})$--$(\ref{ic1})$ 
 with an initial function satisfying
relations $(\ref{le1})$  has the form 
\begin{equation}\label{ue} 
 u(x,t) = \Lambda^{-}_0 \mathrm{erfc} ( \eta)+ 
 \Lambda^{+}_0 \mathrm{erfc} ( -\eta) 
 + \sum\limits_{n=1}^{\infty} t^{-n/2} 
 \left(H_{n,0}(\eta) + H_{n,1}(\eta)\ln
t \right), 
 \end{equation} 
 where $H_{n,m}$ are $C^{\infty}$-smooth functions
of the self-similar variable $$ 
 \eta = \frac{x}{2\sqrt{t}}. 
 $$ 
 } 
 
 Notice that $H_{n,1}$ are the Hermite functions:
 $$ 
 H_{n,1}(\eta) = 
 \frac{\Lambda^+_n - \Lambda^-_n}{4\sqrt{\pi}} 
 \sum\limits_{k=0}^{[(n-1)/2]} 
 \frac{\eta^{n-2k-1} e^{-\eta^2}}{4^k
k! (n-2k-1)!}. 
 $$ 
 
 The above reasoning without
 essential changes can be
applied to the case when $u(x,0)$
grows power-like as $x\to\infty$. 
 For a function $\Lambda$, 
 satisfying the asymptotic relations $$ 
 \Lambda(x) = x^{p}\sum\limits_{n=0}^{\infty} 
 \frac{\Lambda^{\pm}_n}{x^{n}}, 
 \qquad 
 x\to \pm\infty, 
 $$ 
 where~$p$ is a natural number, 
the expansion of the solution has the following form: 
 $$ 
 u(x,t) = 
 \sum\limits_{m=0}^{p} t^{m/2} \left[ 
 \Lambda^{-}_{p-m} \Pi^{-}_{m}(\eta)\, \mathrm{erfc} ( \eta)+ 
 \Lambda^{+}_{p-m} \Pi^{+}_{m}(\eta)\, \mathrm{erfc} ( -\eta)\,+ 
 \right. 
 $$ 
 $$ 
 \left. 
 + P_{m-1}(\eta) \exp(-\eta^2) 
 \right] 
 + \sum\limits_{n=1}^{\infty} t^{-n/2} 
 \left[ G_{n,0}(\eta) + G_{n,1}(\eta)\ln
t \right], 
 $$ 
 where $\Pi^{\pm}_{m}(\eta)$ are polynomials of degree~$m$, 
 whose coefficients are  constants,
and $P_{m-1}(\eta)$ are polynomials of degree $m-1$ 
 ($P_{-1}(\eta)\equiv 0$), 
 $G_{n,k}$ are smooth functions. 
 
\vskip 7mm 
 \setcounter{equation}{0} 
\section{Heat equation on a plane} 
 
 Consider the Cauchy problem for the heat equation on
a plane with a locally Lebesgue-integrable 
 initial function $\Lambda:\mathbb{R}^2 \to \mathbb{R}$ 
of slow growth: 
 \begin{eqnarray} 
 \label{eq} 
 \frac{\partial u}{\partial t} = 
 \frac{\partial^2 u}{\partial x_1^2} 
 +\frac{\partial^2 u}{\partial x_2^2}, 
 & 
 \quad t > 0, \\ 
 \label{ic}\phantom{\frac{1}{1}} 
 u(x_1,x_2,0) = \Lambda (x_1,x_2),& 
 \quad (x_1,x_2)\in\mathbb{R}^2. 
 \end{eqnarray} 
 The solution of this problem has the form 
\begin{equation}\label{us} 
 u(x_1,x_2,t) = \frac{1}{4\pi t} 
 \int\limits_{-\infty}^{+\infty} 
 \int\limits_{-\infty}^{+\infty} 
 \Lambda(s_1,s_2) \exp\left\{ - \frac{(s_1-x_1)^2+(s_2-x_2)^2}{4t} 
 \right\}\, 
 ds_1 ds_2. 
 \end{equation} 
 
 Let us construct a uniform asymptotic expansion
as $|x_1| +|x_2| + t \to \infty$ 
 under the assumption that 
\begin{equation}\label{lmc} 
 \Lambda(x_1,x_2) = 0, \quad x_1<0, 
 \end{equation} 
 \begin{equation}\label{le} 
 \Lambda(x_1,x_2) = x_1^p \sum\limits_{n=0}^{\infty} 
 \frac{\Lambda_n(x_2)}{x_1^{n}}, 
 \quad 
 x_1\to +\infty, 
 \end{equation} 
 where~$p$ is a nonnegative integer, $\Lambda_n$
are continuous functions. 
 In addition, we assume that 
\begin{equation}\label{sc} 
 \begin{split} 
 &\mathrm{supp}\,\Lambda \subset 
 \{ (x_1,x_2): x_1>0,\ |x_2|<|x_1|^{\nu},\ \nu>0\},\\ 
 &\mathrm{supp}\,\Lambda_n\subset [-R_n,R_n],\ R_n>0. 
 \end{split} 
 \end{equation} 
 
 Taking into account condition~(\ref{lmc}), 
we represent function~(\ref{us}) in the form 
\begin{equation}\label{ud} 
 u(x_1,x_2,t)= U_0(x_1,x_2,t) + U_1(x_1,x_2,t), 
 \end{equation} 
 where $$ 
 U_0(x_1,x_2,t) = \int\limits_{0}^{\sigma}
ds_1 
 \int\limits_{-\infty}^{+\infty} ds_2 
 \, \dots \,, 
 \qquad 
 U_1(x_1,x_2,t) = \int\limits_{\sigma}^{+\infty}
ds_1 
 \int\limits_{-\infty}^{+\infty} ds_2 
 \, \dots \,, 
 $$ 
 \begin{equation}\label{svar} 
 \sigma = (x_1^2 + x_2^2+ t)^{\beta/2}, 
 \qquad 
 0 < \beta < 1, 
 \end{equation} 
 and dots denote the integrand from formula~(\ref{us}) 
 together with the factor $(4\pi t)^{-1}$. 
 In the integral $U_1(x_1,x_2,t)$
 we make the change $s_1 = 2z\sqrt{t}$. 
 Setting 
\begin{equation}\label{mevar} 
 \mu= \frac{\sigma}{2\sqrt{t}}, 
 \qquad 
 \eta_1 = \frac{x_1}{2\sqrt{t}} 
 \end{equation} 
 and using by condition~(\ref{le}), we
obtain $$ 
 U_1(x_1,x_2,t) = \frac{1}{2\pi \sqrt{t}} 
 \int\limits_{\mu}^{+\infty} 
 e^{-(\eta_1-z_1)^2} 
 \int\limits_{-\infty}^{+\infty} 
 \Lambda(2z\sqrt{t},s_2) 
 \exp\left\{ - \frac{(s_2-x_2)^2}{4t} \right\}\, ds_2
dz= 
 $$ 
 \begin{equation}\label{ui} 
 = \frac{t^{p/2}}{\sqrt{\pi}} 
 \sum\limits_{n=0}^{N-1} 2^{p-n} t^{-n/2} 
 \int\limits_{\mu}^{+\infty} z^{p-n} 
 e^{-(z-\eta_1)^2} dz 
 \times 
 \end{equation} 
 $$ 
 \times 
 \frac{1}{2\sqrt{\pi t}} 
 \int\limits_{-\infty}^{+\infty} 
 \Lambda_n(s_2) \exp\left\{ - \frac{(s_2-x_2)^2}{4t} \right\}\, ds_2 
 +O( \sigma^{- \rho_1 N}), 
 \qquad \rho_1>0. 
 $$ 
 
 According to Theorem~1, 
from condition~(\ref{sc}) we have 
\begin{equation}\label{s2} 
 \begin{split} 
 \frac{1}{2\sqrt{\pi t}} 
 & \int\limits_{-\infty}^{+\infty} 
 \Lambda_n(s_2) 
 \exp\left\{ - \frac{(s_2-x_2)^2}{4t} \right\}\,
ds_2=\\ 
 & = \exp(-\eta_2^2) 
 \sum\limits_{m=1}^{N}t^{-m/2} 
 Q_{n,m-1}(\eta_2) + 
 O((x_2^2+t)^{-\rho_3 N}), 
 \qquad \rho_3>0, 
 \end{split} 
 \end{equation} 
 where $Q_{n,m-1}(\eta_2)$ are polynomials of degree $m-1$
 in  $\eta_2 = \displaystyle\frac{x_2}{2\sqrt{t}}$, 
 whose coefficients depend on $n$. 
 Then 
\begin{equation}\label{u1} 
 U_1(x_1,x_2,t)= t^{p/2} 
 \sum\limits_{n=1}^{N} t^{-n/2} 
 \widetilde{S}_n(\eta_1,\eta_2) + 
 V_1(\mu,\eta_1,\eta_2,t) + 
 O(\sigma^{-\rho_4 N}), \quad\rho_4>0, 
 \end{equation} 
 \begin{equation}\label{v1f} 
 V_1(\mu,\eta_1,\eta_2,t)= 
 \exp\left\{-(\eta_1^2+\eta_2^2) \right\} 
 \sum\limits_{r^2_s+q^2_s\neq 0} 
 a'_s \eta_1^{m_s}\eta_2^{n_s} t^{l_s} \mu^{r_s} \ln^{q_s}\mu, 
 \end{equation} 
 where $a'_s$ are constants. 
 Coefficients $\widetilde{S}_n(\eta_1,\eta_2)$  are
smooth functions; 
 in particular, $$ 
 \widetilde{S}_1(\eta_1,\eta_2) = 
 \exp\left\{-(\eta_1^2+\eta_2^2) \right\} 
 \left[ 
 \Pi^{(1)}_p(\eta_1)\exp(\eta_1^2)\mathrm{erfc}(-\eta_1) 
 +\Pi_{p-1}^{(2)}(\eta_1) \right], 
 $$ 
 where $\Pi^{(1)}_p(\eta_1)$ and $\Pi_{p-1}^{(2)}(\eta_1)$ 
are polynomials of degree~$p$ and $p-1$,
respectively. 
 
 Let us represent the integral $U_0(x_1,x_2,t)$ 
 in the following form: 
 $$ 
 U_0(x_1,x_2,t) = 
 \frac{\exp\left\{-(\eta_1^2+\eta_2^2) \right\} }{4\pi
t} 
 \times 
 $$ 
 $$ 
 \times 
 \int\limits_{0}^{\sigma} 
 \int\limits_{-\infty}^{\infty} 
 \Lambda(s_1,s_2) 
 \sum\limits_{m=0}^{N} \frac{1}{m!} 
 \left( \frac{\eta_1 s_1+\eta_2 s_2}{\sqrt{t}} - \frac{s_1^2+s_2^2}{4t}\right)^m 
 ds_2 ds_1 + 
 O( \sigma^{- \rho_5 N}), 
 \qquad \rho_5 > 0. 
 $$ 
 From this formula we obtain $$ 
 U_0(x_1,x_2,t)= \exp\left\{-(\eta_1^2+\eta_2^2) \right\} 
 \sum\limits_{n=2}^{N} t^{-n/2} 
 \sum\limits_{\substack{ 0\leqslant m_1+m_2\leqslant
n-2\\ 0\leqslant l_1+l_2\leqslant
n-2}} 
 a_{m_1,m_2,l_1,l_2} \eta_1^{m_1}\eta_2^{m_2} 
 \times 
 $$ 
 $$ 
 \times 
 \int\limits_{0}^{\sigma} 
 \int\limits_{-\infty}^{\infty} 
 s_1^{l_1} s_2^{l_2} \Lambda(s_1,s_2) ds_2
ds_1 
 +O( \sigma^{- \rho_5 N}), 
 $$ 
 where $a_{m_1,m_2,l_1,l_2}$ are some
constants. 
 Let us transform the integral as follows: 
 $$ 
 \int\limits_{0}^{\sigma} 
 \int\limits_{-\infty}^{\infty} 
 s_1^{l_1} s_2^{l_2} \Lambda(s_1,s_2) ds_2
ds_1 = 
 \int\limits_{0}^{1} 
 \int\limits_{-\infty}^{\infty} 
 s_1^{l_1} s_2^{l_2} \Lambda(s_1,s_2) ds_2
ds_1 + 
 $$ 
 $$ 
 +\int\limits_{1}^{\sigma} 
 \int\limits_{-\infty}^{\infty} 
 s_1^{l_1} s_2^{l_2}\left[ \Lambda(s_1,s_2) 
 -s_1^p\Lambda_0(s_2) - \dotsc - s_1^{-l_1-1} \Lambda_{p+l_1+1}(s_2) 
 \right] 
 ds_2 ds_1 + 
 $$ 
 $$ 
 +\int\limits_{1}^{\sigma} 
 \int\limits_{-\infty}^{\infty} 
 \left[ 
 s_1^{l_1+p} s_2^{l_2}\Lambda_0(s_2) + \dotsc +
s_1^{-1}s_2^{l_2} \Lambda_{p+l_1+1}(s_2) 
 \right] 
 ds_2 ds_1 = 
 $$ 
 $$ 
 = A_{l_1,l_2} + \ln\sigma 
 \int\limits_{-\infty}^{\infty} 
 s_2^{l_2} \Lambda_{p+l_1+1}(s_2) ds_2 
 + \sum\limits_{k=1}^{N-1}(c_{k}\sigma^{k}+c_{-k}\sigma^{-k}) 
 + O(\sigma^{-N}) = 
 $$ 
 $$ 
 = A_{l_1,l_2} + B_{l_1,l_2} \ln t 
 +2B_{l_1,l_2}\ln(2\mu) 
 + \sum\limits_{k=1}^{N-1}(c'_{k}\mu^{k}t^{k/2}+c'_{-k}\mu^{-k}t^{-k/2}) 
 + O(\sigma^{-N}), 
 $$ 
 where $A_{l_1,l_2}$, $B_{l_1,l_2}$, $c'_{k}$
and $c'_{-k}$ are constants. 
 Then $$ 
 U_0(x_1,x_2,t) = 
 \exp\left\{-(\eta_1^2+\eta_2^2) \right\} 
 \sum\limits_{n=2}^{N} t^{-n/2} 
 \left[ \Pi_{n-2}(\eta_1,\eta_2) + 
 \Pi_{n-2}^{*}(\eta_1,\eta_2)\ln t 
 \right] + 
 $$ 
 \begin{equation}\label{u0} 
 + V_0(\mu,\eta_1,\eta_2,t) 
 +O( \sigma^{- \rho_6 N}), 
 \qquad \rho_6>0, 
 \end{equation} 
 where $\Pi_{n-2}(\eta_1,\eta_2)$ and $\Pi_{n-2}^{*}(\eta_1,\eta_2)$ 
are polynomials of degree $n-2$, 
 and the expression $$ 
 V_0(\mu,\eta_1,\eta_2,t)= 
 \exp(-\eta_1^2-\eta_2^2) 
 \sum\limits_{r^2_s+q^2_s\neq 0} 
 a''_s \eta_1^{m_s}\eta_2^{n_s} t^{l_s} \mu^{r_s} \ln^{q_s}\mu, 
 $$ 
 where $a''_s$ are some constants,
 is obtained similarly to expression~(\ref{v1f}). 
 Expressions $V_{0}(\mu,\eta,t)$ and $V_{1}(\mu,\eta,t)$ 
 asymptotically diminish 
$$ 
 V_{0}(\mu,\eta_1,\eta_2,t)+V_{1}(\mu,\eta_1,\eta_2,t) 
 =O( \sigma^{- \rho N}), 
 \qquad 
 (x_1,x_2,t)\in T_{\alpha}, 
 $$ 
 where $\rho=\min\{ \rho_4,\rho_6 \}$. 
 Then we arrive at the following statement 
(details of the proof see in~\cite{hd2}).

\vskip 5mm 
 \textbf{Theorem~2.} 
 {\it 
 As $|x_1| + |x_2| + t \to \infty$ the asymptotics 
of the solution of equation $(\ref{eq})$ 
 with conditions $(\ref{ic}),$ $(\ref{lmc}),$ $(\ref{le})$ and $(\ref{sc})$ 
 has the form $$ 
 u(x_1,x_2,t) = 
 \sum\limits_{n=1}^{\infty} t^{-n/2} 
 \left[ t^{p/2} S_{n}(\eta_1,\eta_2) + 
 t^{-1/2}\ln t\, \Pi_{n}(\eta_1,\eta_2) \exp(-\eta_1^2-\eta_2^2) 
 \right], 
 $$ 
 where $S_{n}(\eta_1,\eta_2)$ are
 $C^{\infty}$-smooth functions 
of slow growth,
 $\Pi_{n}(\eta_1,\eta_2)$ are polynomials of degree~$n$ 
 in the self-similar variables $$ 
 \eta_1 = \frac{x_1}{2\sqrt{t}}, 
 \qquad 
 \eta_2 = \frac{x_2}{2\sqrt{t}}. 
 $$ 
 } 
 \bigskip 
 
\vskip 7mm
 \setcounter{equation}{0} 
 \section{Multidimensional equation}
 
 For a multidimensional problem,
 there holds the following result. 
 
 \vskip 5mm
 \textbf{Theorem~3.} 
 {\it Let $u(x_1,\dotsc,x_m,t)$ be the solution
of the Cauchy problem
\begin{eqnarray*} 
 \frac{\partial u}{\partial t} = 
 \frac{\partial^2 u}{\partial x_1^2} 
 +\dotsc 
 +\frac{\partial^2 u}{\partial x_m^2}, 
 & 
 \quad t > 0, \\ 
 \phantom{\frac{1}{1}} 
 u(x_1,\dotsc,x_m,0) = \Lambda (x_1,\dotsc,x_m),& 
 \quad (x_1,\dotsc,x_m)\in\mathbb{R}^m, 
 \end{eqnarray*} 
 with a locally integrable initial function $\Lambda$ 
of slow growth. 
 If the conditions $$ 
 \Lambda(x_1,\dotsc,x_m) = 0, \quad
x_1<0, 
 $$ 
 $$ 
 \Lambda(x_1,\dotsc,x_m) = x_1^p \sum\limits_{n=0}^{\infty} 
 \frac{\Lambda_n(x_2,\dotsc,x_m)}{x_1^{n}}, 
 \quad 
 x_1\to +\infty, 
 $$ 
 where $$ 
 \mathrm{supp}\,\Lambda \subset 
 \{ (x_1,\dotsc,x_m): x_1>0,\ |x_2| + \dotsc + |x_m|<|x_1|^{\nu},\ \nu>0\}, 
 $$ 
 $$ 
 \mathrm{supp}\,\Lambda_n\subset [-R_n,R_n]^{m-1},\quad R_n>0, 
 $$ 
are fulfilled, then
 there holds the asymptotic formula $$ 
 u(x_1,\dotsc,x_m,t) = t^{-m/2} 
 \sum\limits_{n=0}^{\infty} t^{-n/2} 
 \left[ t^{(p+1)/2} S_{n}(\eta_1,\dotsc,\eta_m) + 
 \right. 
 $$ 
 $$ 
 \phantom{t^{1/2}} 
 \left.+\ln t\, \Pi_{n}(\eta_1,\dotsc,\eta_m) 
 \exp(-\eta_1^2-\dotsc-\eta_m^2) 
 \right], 
 \qquad 
 |x_1| + \dotsc + |x_m| + t \to \infty, 
 $$ 
 where $S_{n}(\eta_1,\dotsc,\eta_m)$ are smooth
functions of slow growth,
 $\Pi_{n}(\eta_1,\dotsc,\eta_m)$ are polynomials of degree~$n$ 
 in the self-similar variables $$ 
 \eta_1 = \frac{x_1}{2\sqrt{t}}, 
 \quad 
 \dotsc, 
 \quad 
 \eta_m = \frac{x_m}{2\sqrt{t}}\, . 
 $$ 
 }

\end{document}